\documentclass[12pt]{article}
\def\today{January 6, 2004}
\usepackage{amsmath}
\usepackage{amssymb}
\usepackage{amsthm}
\usepackage[pdfpagemode=None]{hyperref}
\usepackage{psfig}
\usepackage{graphicx}
\newtheorem{thm}{Theorem}[section]
\newtheorem{co}[thm]{Corollary}
\newtheorem{lem}[thm]{Lemma}

\newtheorem{assumption}[thm]{Assumption}

\newtheorem{definition}[thm]{Definition}
\newenvironment{de}{\begin{definition}\rm}{\end{definition}}
\newtheorem{example}[thm]{Example}

\newtheorem{remark}[thm]{Remark}
\newenvironment{rem}{\begin{remark}\rm}{\end{remark}}
\newtheorem{tab}{Table}
\newenvironment{ta}{\begin{tab}\rm}{\end{tab}}

\newcommand{\Section}[1]{\section{#1}\setcounter{equation}{0}}

\newcommand{\eqr}[1]{~\mbox{$(${\rm \ref{#1}}$)$}}
\newcommand{\tr}{{\rm tr}\,}
\newcommand{\diag}{{\rm diag}\,}
\newcommand{\tang}{{\rm tang}\,}

\newcommand{\dist}{{\rm dist}\,}

\newcommand{\V}{\mathcal{V}}

\title{Unitary Space Time Constellation Analysis: \\
An Upper Bound for the Diversity \footnote{Both authors were
supported in part by NSF grants
    DMS-00-72383 and CCR-02-05310. The first author was also
    supported by a fellowship from the Center of Applied
    Mathematics at the University of Notre Dame.}}

\date{{\normalsize \today}}

\author{Guangyue Han,\ \ Joachim Rosenthal\\
  {\normalsize Department of Mathematics}\vspace{-1mm} \\
  {\normalsize University of Notre Dame}\vspace{-1mm} \\
  {\normalsize Notre Dame, IN 46556.}\\
  {\normalsize {\em e-mail:\/} Han.13@nd.edu,\
  Rosenthal.1@nd.edu}\vspace{-1mm}\\
  {\normalsize {\em URL:} http://www.nd.edu/\~{}eecoding/} }

\begin{document}\maketitle\thispagestyle{empty}
\begin{abstract}
  The diversity product and the diversity sum are two very
  important parameters for a good-performing unitary space time
  constellation. A basic question is what the maximal diversity
  product (or sum) is. In this paper we are going to derive
  general upper bounds on the diversity sum and the diversity
  product for unitary constellations of any dimension $n$ and any
  size $m$ using packing techniques on the compact Lie group
  $U(n)$.
\end{abstract}
\Section{Introduction}

Let $A$ be a matrix with complex entries. $A^*$ denotes the
conjugate transpose of $A$. Let $\| \hspace{2mm} \|$ denote the
Frobenius norm of a matrix, i.e.,
$$
\|A\|=\sqrt{\tr(AA^*)}.
$$
A square matrix $A$ is called unitary if $A^*A=AA^*=I$, where $I$
denotes the identity matrix. We denote by $U(n)$ the set of all
$n \times n$ unitary matrices. $U(n)$ is a real algebraic variety
and a smooth manifold of real dimension $n^2$. For the purpose of
this paper a unitary space time constellation (or code) $\V$ is
simply a finite subset of $U(n)$,
$$
\V=\{A_1, A_2, \cdots, A_m\} \subset U(n).
$$
We say $\V$ has dimension $n$ and size $m$. Unitary space time
codes have been intensely studied in recent years and we refer
the interested readers to~\cite{al98,ho00,ho00a,sh01} and the
references of these papers. The readers will find the motivation
and engineering applications of such kind of codes. The quality
of a unitary space time code is governed by two important
parameters, the diversity product and the diversity sum.

\begin{de}
  The {\em diversity product}~\cite{ho00} of a unitary space time
  code $\V$ is defined through
  $$
  \prod \V:=\frac{1}{2} \min \{{|\det(A-B)|}^{\frac{1}{n}} |A,
  B \in \V, A \neq B\}.
  $$
  The {\em diversity sum}~\cite{li02} is defined as
  $$
  \sum \V:=\frac{1}{2\sqrt{n}} \min \{\|A-B\| |A, B \in \V, A
  \neq B \}.
  $$
\end{de}
$\V$ is called fully diverse if $\prod \V > 0$. As explained
in~\cite{ha03u4}, a space time code with large diversity sum
tends to perform well at low signal to noise ratios whereas a
code with a large diversity product tends to perform well at high
signal to noise ratios. A major coding design problem is the
construction of unitary space time codes where the diversity sum
(or product) is optimal or near optimal inside the set of all the
space time codes with the same parameters $n, m$. We would like
to remark that for every positive integer $n$ and $m$, a Haar
distributed random space time code is fully diverse with
probability $1$. A simple proof can be found in~\cite{ha03u4}.

The purpose of this paper is to derive for $n$ and $m$ tight
upper bounds for the diversity product $\prod \V$ and the
diversity sum $\sum \V$. When $n=1$ then trivially
$|\det(A-B)|=\| A-B\|$ and it follows that $\sum \V=\prod \V$ in
this situation.  The following lemma states that for every space
time code $\V$, $\sum \V$ is an upper bound for $\prod \V$ and by
having an upper bound for $\sum \V$ we immediately also have an
upper bound for $\prod \V$. The readers can find the statements
about the relationship between $\prod \V$ and $\sum \V$
in~\cite{li02}, for completeness we include a detailed proof.

\begin{lem}
  For any unitary space time code $\V$,
  $$
  \prod \V \leq \sum \V.
  $$
\end{lem}

\begin{proof}
  Let $C$ be an $n \times n$ complex matrix with singular value
  decomposition
  $$
  C=U \diag(c_1, c_2, \cdots, c_n) V,
  $$
  where $U, V$ are unitary matrices and $c_j \geq 0$ for $j=1,
  2, \cdots, n$ are the singular values of $C$. First we are going to prove
  $$
  \frac{1}{2} {|\det(C)|}^{\frac{1}{n}} \leq \frac{1}{2
    \sqrt{n}} \|C\|.
  $$
  If $c_j=0$ for some $j$, then the inequality is trivial. Hence we assume $c_j
  >0$ for all $j$'s. Because $U, V$ are unitary matrices, it
  follows that
  $$
  \frac{1}{2} {|\det(C)|}^{\frac{1}{n}}=\frac{1}{2}
  \left(\prod_{j=1}^{n} c_j \right)^{\frac{1}{n}}.
  $$
  Similarly one verifies that
  $$
  \frac{1}{2 \sqrt{n}} \|C\|= \frac{1}{2 \sqrt{n}}
  \sqrt{\sum_{j=1}^n c_j^2}.
  $$
  Applying Cauchy-Schwarz inequality, we have
  $$
  \left(\prod_{j=1}^{n} c_j \right)^{\frac{1}{n}} \leq
  \frac{\sum_{j=1}^n c_j}{n} \leq \frac{1}{\sqrt{n}}
  \sqrt{\sum_{j=1}^n c_j^2}.
  $$
  Hence one concludes that for an $n \times n$ square matrix
  $C$,
  $$
  \frac{1}{2} {|\det(C)|}^{\frac{1}{n}} \leq \frac{1}{2
    \sqrt{n}} \|C\|.
  $$
  By the definition of $\prod \V$, $\sum \V$ and the above
  inequality one gets
  $$
  \prod \V \leq \sum \V.
  $$
\end{proof}

Of course it would be desirable to know for every $n$ and $m$ what
the largest possible value of $\sum \V$ is. This is the
motivation of the following definition.

\begin{de}
  Let $\Delta(n, m)$ be the infimum of all numbers such that for
  every unitary space time code $\V$ of dimension $n$ and size
  $m$, one has
  $$
  \sum \V \leq \Delta (n, m).
  $$
\end{de}
\begin{rem}
  As pointed out by Liang and Xia~\cite{li02} there exists a
  constellation $\V$ of dimension $n$ and size $m$ with $\sum
  \V=\Delta(n,m)$. This is due to the fact that $U(n)^m$ is a
  compact manifold.
\end{rem}

The exact values of $\Delta (n, m)$ are only known in very few
special cases. In the case $n=1$, one checks that $\Delta (1,
m)=\sin \frac{\pi}{m}$ for $m \geq 2$. When $n \geq 2$ and $m=3$,
one has $\Delta (n, 3)=\frac{\sqrt{3}}{2}$. When $m=2$, we have
$\Delta (n, 2)=1$ for $n \geq 2$. For $n=2$, the following values
were computed in~\cite{li02}.
\begin{center}
\begin{tabular}{|c|c|c|c|c|c|c|c|c|c|} \hline
  $m$ & $2$ & $3$ & $4$ & $5$ & $6$ & $7$ & $8$ & $9$ & $10$ through $16$ \\
  \hline
  $\Delta  (2, m)$ & $1$ & $\frac{1}{2}\sqrt{3}$ &
  $\frac{1}{3}\sqrt{6}$ & $\frac{1}{4}\sqrt{10}$ &
  $\frac{1}{5}\sqrt{15}$ & $\frac{1}{6}\sqrt{21}$ &
  $\frac{1}{7}\sqrt{28}$ & $\frac{1}{8}\sqrt{36}$ & $\frac{1}{2}\sqrt{2}$ \\
  \hline
\end{tabular}
\end{center}

Liang and Xia~\cite{li02} observed the connection between a
unitary constellation and an Euclidean sphere code and
beautifully derived an upper bound for $2$ dimensional unitary
constellations which is very tight when $m \leq 100$.  In this
paper we present a new general upper bound for $\Delta (n, m)$
for every dimension $n$ and every size $m$ while improving
certain results in~\cite{li02}. To the best of our knowledge the
new upper bounds we derived are tighter than any previously
published bounds as soon as $m$ is sufficiently large.

\Section{Upper Bound Analysis} \label{Upper-Bound}

In this section we are going to study the packing problem on
$U(n)$ and derive three upper bounds for the numbers $\Delta
(n, m)$. All the resulted bounds are derived by differential
geometric means and all bounds can be viewed as certain sphere
packing bounds.

From a differential geometry point of view we can view $U(n)$ as
a $n^2$-dimensional compact Lie group. $U(n)$ is also naturally a
submanifold of the Euclidean space $\mathbb{R}^{2n^2}$. In this
way $U(n)$ will have the induced geometry of the standard
Euclidean geometry of $\mathbb{R}^{2n^2}$. Finally there is a
third way to see $U(n)$ as a submanifold of another Riemannian
manifold $S(n)$ and we will say more later.

The basic strategy for computing the upper bounds for $\Delta (n,
m)$ is as follows. Given a unitary space time code $\V=\{A_1,
A_2, \cdots, A_m\}$, around each matrix $A_j$ we can choose a
neighborhood $N_r(A_j)$ with radius $r$ (the radius will be
specified later). Let $V_j=V(N_r(A_j))$ be the volume of the
neighborhood $N_r(A_j)$. If all the neighborhoods are
non-overlapping, then necessarily we will have
$$
\sum_{j=1}^m V_j \leq V(U(n)),
$$
where $V(U(n))$ denotes the total volume of unitary group $U(n)$.
This inequality in turn will result in an upper bound for the
numbers $\Delta  (n, m)$. By employing different metrics
(Euclidean or Riemannian) and by considering different embeddings
of $U(n)$, we derive three different upper bounds for $\Delta (n,
m)$.

Let $\mathcal{M}_1$ be the manifold consisting of all the $n
\times n$ Hermitian matrices, i.e.
$$
\mathcal{M}_1=\{H| H=H^*\}.
$$
$\mathcal{M}_1$ has dimension $n^2$ and can be viewed
isometrically as Euclidean space $\mathbb{R}^{n^2}$. Assume that
$H=(H_{jk})$ and assume that $H_{jk}=x_{jk}+i y_{jk}$. We use
$(dH)$ to denote the volume element of $\mathcal{M}_1$, where
\begin{equation}   \label{Hermitian-volume}
(dH)=\left( \frac{i}{2} \right)^{n(n-1)/2}\bigwedge_{l=1}^n dH_{ll}
\bigwedge_{j < k} dH_{jk} \bigwedge_{j < k}
d\bar{H}_{jk}=\bigwedge_{l=1}^n dx_{ll} \bigwedge_{j < k} dx_{jk}
\bigwedge_{j < k} dy_{jk}.
\end{equation}
With a small abuse of the notation, one can check that the volume
element of $\mathcal{M}_2$, the manifold consisting of all the $n
\times n$ skew-Hermitian matrices, can be written as
\begin{equation}  \label{skew-Hermitian-volume}
(dH)=\left( \frac{i}{2} \right)^{n(n-1)/2} \left( \frac{1}{i}
\right)^n \bigwedge_{l=1}^n dH_{ll} \bigwedge_{j < k} dH_{jk} \bigwedge_{j
< k} d\bar{H}_{jk}=\bigwedge_{l=1}^n dy_{ll} \bigwedge_{j < k} dx_{jk}
\bigwedge_{j < k} dy_{jk}.
\end{equation}

For a unitary matrix $U$, if we differentiate $U^*U=I$, we will
have
$$
U^* dU + dU^* U=0.
$$
Therefore $U^*dU$ is skew-Hermitian. The following lemma will
characterize the volume element of $U(n)$. For the terminologies
in this lemma, we refer to a standard differential geometry or
integral geometry book, e.g. ~\cite{he01, sa76}.

\begin{lem}
  The volume element of $U(n)$ induced by the Euclidean space
  $\mathbb{R}^{2n^2}$ is bi-invariant and the volume element can
  be written as $(U^*dU)$ up to a scalar constant.
\end{lem}

\begin{proof}
  The bi-invariance comes from the orthonormality of $U(n)$.
  $(U^*dU)$ is left-invariant according to the definition. Indeed
  for a fixed yet arbitrary unitary matrix $V$,
  $$
  (V U)^* d(V U)= U^* V^* V dU=U^*dU.
  $$
  Since $U(n)$ is a compact Lie group and any compact Lie
  group is unimodular, $(U^*dU)$ is also right-invariant.
  Because the bi-invariant $n^2$ differential forms are unique up
  to a scalar, one concludes that the volume element can be
  written as $(U^*dU)$.
\end{proof}

The following theorem will represent the volume element of $U(n)$
in another way. One will see that it is closely related to the
eigenvalues of unitary matrices.

\begin{thm}
  For the Schur decomposition of a unitary matrix $\Theta$:
\begin{equation}  \label{Schur}
\Theta=U\diag(e^{i\theta_1},e^{i\theta_2},\cdots,e^{i\theta_n})U^*,
\end{equation}
we will have
\begin{equation}  \label{Schur-volume-element}
(\Theta^* d\Theta)=\prod_{j<k}{|e^{i \theta_j}-e^{i \theta_k}|}^2
d\theta_1 \wedge d\theta_2 \wedge \cdots \wedge d\theta_n \wedge
(U^*dU-\diag(U^*dU)).
\end{equation}
\end{thm}

\begin{proof}
  Let $D=\diag(e^{i\theta_1},e^{i\theta_2},\cdots,e^{i\theta_n})$
  and take the differential of Equation\eqr{Schur},
  $$
  d\Theta=dU D U^*+U dD U^*+ U D dU^*.
  $$
  It follows that,
  $$
  \Theta^*d\Theta=U D^* U^* dU D U^* +U D^* dD U^*+U
  dU^*=U(D^*U^* dU D+ D^*d D) U^*+U dU^*.
  $$
  Due to the right-invariance of the volume element in $U(n)$,
  it follows that
  $$
  (\Theta^*d\Theta)=(U^*\Theta^* d\Theta U)=(D^* U^*dU D -
  U^*dU+i \diag(d\theta_1, d\theta_2, \cdots, d\theta_n)).
  $$
  Note that $(D^* U^*dU D - U^*dU)_{jk}=(e^{i \theta_j}-e^{i
    \theta_k}) U_{jk}$, therefore the diagonal elements of $D^*
  U^*dU D - U^*dU$ are all zeros and the off diagonal elements
  are scaled version of the ones of $U^* dU$.  According to
  formula\eqr{skew-Hermitian-volume}, the claim in the theorem
  follows.
\end{proof}

The following theorem calculates the volume of a small
neighborhood with Euclidean distance $r$. Because of the
homogeneity of $U(n)$, the center of this small ``ball'' is
chosen to be $I$ without loss of generality. For a unitary matrix
$U$, we assume $e^{i \theta_j}$'s are its eigenvalues, i.e., $U
\sim \diag(e^{i \theta_1}, e^{i \theta_2}, \cdots, e^{i
  \theta_n})$. For a fixed unitary matrix $A$, let
$$
U_r^E(n, A)=\{U \in U(n)| \|U-A\| \leq r\}.
$$
Again because of the homogeneity of $U(n)$, $V(U_r^E(n, A))$ does
not depend on the choice of $A$. In the sequel $V(U_r^E(n))$ will
be used to denote $V(U_r^E(n, A))$ for any unitary matrix $A$.
Let $S(n)$ denote a $2n^2-1$ dimensional sphere centered at the
origin with radius $\sqrt{n}$, i.e.,
$$
S(n)=\{(x_1,x_2,\cdots,x_{2n^2})|x_1^2+x_2^2+\cdots+x_{2n^2}^2=n\}.
$$
Apparently $U(n)$ is a submanifold of $S(n)$. For a particular
point $S_0 \in S(n)$, let
$$
S_r(n,S_0)=\{S \in S(n) |\|S-S_0\| \leq r\}.
$$

\begin{thm}  \label{Euclidean-Ball}
  Let
  \begin{equation}              \label{D1}
  D_1=\{(\theta_1, \theta_2, \cdots, \theta_n)| -\pi \leq
  \theta_j < \pi \; \mbox{for} \; j=1, 2, \cdots, n \}
  \end{equation}
  and
  \begin{equation}              \label{D2}
  D_2=\left\{(\theta_1, \theta_2, \cdots, \theta_n)| \sum_{j=1}^n
  \sin^2 \frac{\theta_j}{2} \leq \frac{r^2}{4} \right\},
   \end{equation}
  then
\begin{equation}   \label{Euclidean-Volume}
V(U_r^E(n))=\frac{\int\!\!\!\int_{D_1 \cap D_2} \prod_{j<k}{|e^{i
\theta_j}-e^{i \theta_k}|}^2 d\theta_1 d\theta_2 \cdots d\theta_n
}{\int\!\!\!\int_{D_1} \prod_{j<k}{|e^{i \theta_j}-e^{i
\theta_k}|}^2 d\theta_1 d\theta_2 \cdots d\theta_n} V(U(n)).
\end{equation}
\end{thm}

\begin{proof}
  Note that ${\|I-U\|}_2 \leq r$ is equivalent to $\sum_{j=1}^n
  \sin^2 \frac{\theta_j}{2} \leq \frac{r^2}{4}$. For a given
  unitary matrix $\Theta$, the Schur decomposition $\Theta=U^*
  \diag(e^{i \theta_1}, e^{i \theta_2}, \cdots, e^{i \theta_n})
  U$ is unique if $\theta_j$'s are strictly ordered. So if we
  take the integral of formula\eqr{Schur-volume-element} over the
  integration region disregarding the order of $\theta_j$'s, we
  will obtain $n!$ times the volume of $V(U_r^E(n))$. Thus the
  volume of $U_r^E(n)$ will be
  $$
  V(U_r^E(n))=\frac{1}{n!} \int\!\!\!\int_{D_1 \cap D_2}
  \prod_{j<k}{|e^{i \theta_j}-e^{i \theta_k}|}^2 d\theta_1
  d\theta_2 \cdots d\theta_n \int\!\!\!\int_{U(n)}
  (U^*dU-\diag(U^*dU)).
  $$

  Using the same argument, we will derive the volume of $U(n)$:
  $$
  V(U(n))=\frac{1}{n!} \int\!\!\!\int_{D_1} \prod_{j<k}{|e^{i
      \theta_j}-e^{i \theta_k}|}^2 d\theta_1 d\theta_2 \cdots
  d\theta_n \int\!\!\!\int_{U(n)} (U^*dU-\diag(U^*dU)).
  $$
  Compare the two derived volume formula, the claim in the
  theorem follows.
\end{proof}

\begin{rem}
  By the Weyl denominator formula~\cite{go98} one can replace
  $$
  \int\!\!\!\int_{D_1} \prod_{j<k}{|e^{i \theta_j}-e^{i
      \theta_k}|}^2 d\theta_1 d\theta_2 \cdots d\theta_n
  $$

  with $(2 \pi)^n n!$. We keep it as it is to make the formula
  literally understandable.
\end{rem}

There are several approaches to derive upper bounds for the
diversity sum. The first approach considers $U(n)$ as a
submanifold of $S(n)$, then chooses the non-overlapping
neighborhoods to be small balls with radius $r$ (with regard to
the Euclidean distance).  This will result in the first upper
bound (B1) which we derive in this paper.

\begin{thm} \label{Sphere-Upper-Bound}
  Let $D_1$ and $D_2$ be defined as in\eqr{D1} and \eqr{D2}.
  Assume $r_0^E=r_0^E (n, m)$ is the solution to the following
  equation (with variable $r$):
\begin{equation}   \label{Euclidean-Equality}
m\int\!\!\!\int_{D_1 \cap D_2} \prod_{j<k}{|e^{i \theta_j}-e^{i
\theta_k}|}^2 d\theta_1 d\theta_2 \cdots d\theta_n =
\int\!\!\!\int_{D_1} \prod_{j<k}{|e^{i \theta_j}-e^{i
\theta_k}|}^2 d\theta_1 d\theta_2 \cdots d\theta_n,
\end{equation}
then
\begin{equation}
\Delta  (n, m)\leq \sqrt{\frac{(r_0^E)^2}{n}-\frac{(r_0^E)^4}{4n^2}}. \tag{B1}
\end{equation}
\end{thm}

\begin{proof}
  For a fixed yet arbitrary unitary constellation $\V=\{A_1, A_2,
  \cdots, A_m\}$, consider $m$ small non-overlapping
  neighborhoods $S_r(n, A_j)$ in $S(n)$. We can increase $r$ such
  that there exist $l,k$ such that $S_r(n, A_l)$ and $S_r(n, A_k)$ are
  tangent to each other.  Apparently
  $$
  U_r^E(n, A_j)=S_r(n, A_j) \cap U(n),
  $$
  for any $j$. Since $S_r(n, A_j)$'s are non-overlapping, we
  conclude that $U_r^E(n, A_j)$'s are non-overlapping.  Therefore
  we have
  $$
  \sum_{j=1}^m V(U_r^E(n, A_j)) \leq V(U(n)),
  $$
  that is
  $$
  m V(U_r^E(n)) \leq V(U(n)).
  $$

  One can check that $V(U_r^E(n))$ is an increasing function of
  $r$, so any $r$ satisfying the above inequality will be less
  than the solution to the equality:
  $$
  m V(U_r^E(n)) = V(U(n)),
  $$
  which is essentially Equality\eqr{Euclidean-Equality}. So we
  conclude that $r \leq r_0^E$.

  Note that any two points $S_0, S_1 \in S(n)$ with two
  non-overlapping neighborhoods $S_r(n, S_0)$ and $S_r(n, S_1)$
  will have distance $\|S_0-S_1\| \geq 2\sqrt{r^2-r^4/(4n)}$,
  where the equality holds only if $S_r(n, S_0)$ and $S_r(n,
  S_1)$ are tangent to each other. Apply the argument to $A_j$'s
  and note that $A_l$ and $A_k$ are the closest pair of points
  with $\|A_l-A_k\|=2\sqrt{r^2-r^4/(4n)}$, we reach the
  conclusion of the theorem.
\end{proof}

For a fixed $S_0 \in S(n)$, consider $S_r(n, S_0) \subset S(n)$.
Let $\tau=\tau(n,r)$ denote the maximal number $\tau$ such that
$S_r(n, S_1),S_r(n, S_2),\cdots,S_r(n, S_{\tau})$ are
non-overlapping and $S_r(n,S_j)$ is tangent to $S_r(n,S_0)$ for
$j=1, 2, \cdots, n$. One checks that $\tau(n, r)$ does not depend
on the choice of $S_0$. In this sense $\tau(n, r)$ can be viewed
as generalized kissing number~\cite{co93b} on an Euclidean sphere.
For a fixed $n$ dimensional unitary constellation $\V=\{A_1, A_2,
\cdots, A_m\}$, let $r(\V)$ denote the maximal radius $r$ such
that $S_r(n, A_1),S_r(n, A_2),\cdots,S_r(n, A_m)$ are
non-overlapping. Let $r_{opt}=r_{opt}(n,m)$ denote the maximal
$r(\V)$ over all possible $n$ dimensional unitary constellation
$\V$ with cardinality $m$. One checks $\Delta
(n,m)=r_{opt}(n,m)/\sqrt{2n}$. The following theorem and
corollary give a lower bound for the optimal diversity sum
$\Delta  (n, m)$.

\begin{thm}
  Let $D_1$  be defined as in\eqr{D1} and
assume that $r_0^E=r_0^E(n, m)$ is the solution to the
  equation\eqr{Euclidean-Equality}. Let
  $$
  \tilde{D}_2=\left\{(\theta_1, \theta_2, \cdots, \theta_n)| \sum_{j=1}^n
  \sin^2 \frac{\theta_j}{2} \leq \frac{(r_0^E)^2}{4} \right\}
  $$
  and let
  $$
  D_3=\left\{(\theta_1, \theta_2, \cdots, \theta_n)|
    \sum_{j=1}^n \sin^2 \frac{\theta_j}{2} \leq r_{opt}(n,
    m)^2-r_{opt}(n,m)^4/(4n) \right\}.
  $$
  Then
  \begin{multline*}
    \int\!\!\!\int_{D_1 \cap \tilde{D}_2} \prod_{j<k}{|e^{i
        \theta_j}-e^{i \theta_k}|}^2 d\theta_1 d\theta_2 \cdots
    d\theta_n  \\
    \leq (\tau(n, r_{opt}(n, m))+1)\int\!\!\!\int_{D_1 \cap D_3}
    \prod_{j<k}{|e^{i \theta_j}-e^{i \theta_k}|}^2 d\theta_1
    d\theta_2 \cdots d\theta_n.
  \end{multline*}
\end{thm}

\begin{proof}
  According to the derivation of $r_0^E$, we have
\begin{equation}  \label{formula-1}
m \int\!\!\!\int_{D_1 \cap \tilde{D}_2} \prod_{j<k}{|e^{i \theta_j}-e^{i
\theta_k}|}^2 d\theta_1 d\theta_2 \cdots
d\theta_n={\int\!\!\!\int_{D_1} \prod_{j<k}{|e^{i \theta_j}-e^{i
\theta_k}|}^2 d\theta_1 d\theta_2 \cdots d\theta_n}.
\end{equation}

Assume that $\V=\{A_1, A_2, \cdots, A_m\}$ is an $n$ dimensional
unitary constellation reaching $r_{opt}(n, m)$, i.e.,
$r(\V)=r_{opt}(n, m)$. For simplicity let $r=r(\V)$. Let $m'$
denote the maximal number such that $S_r(n, A_1), S_r(n, A_2),
\cdots, S_r(n, A_m), \cdots, S_r(n, A_{m'})$ are non-overlapping.
Let $r_1=2\sqrt{r^2-r^2/(4n)}$, we claim that
$$
U(n) \subset \bigcup_{j=1}^{m'} U_{r_1}^E(n, A_j).
$$
Otherwise suppose there is a unitary matrix $A_0 \notin
\bigcup_{j=1}^{m'} U_{r_1}^E(n, A_j)$, then $\|A_0-A_j\| > r_1$
(see Theorem~\ref{Sphere-Upper-Bound}). Thus $S_r(n, A_0)$ does
not intersect with $S_r(n, A_j)$ for $j=1, 2, \cdots, m'$.
Therefore one can find $m'+1$ small balls with radius $r$ which
are non-overlapping. This contradicts the maximality of $m'$.
Thus we have $\sum_{j=1}^{m'} V(U_{r_1}^E(n, A_j)) \geq V(U(n))$,
that is
\begin{equation} \label{formula-2}
m'\int\!\!\!\int_{D_1 \cap D_3} \prod_{j<k}{|e^{i \theta_j}-e^{i
\theta_k}|}^2 d\theta_1 d\theta_2 \cdots d\theta_n \geq
{\int\!\!\!\int_{D_1} \prod_{j<k}{|e^{i \theta_j}-e^{i
\theta_k}|}^2 d\theta_1 d\theta_2 \cdots d\theta_n}.
\end{equation}

We further claim that
\begin{equation} \label{formula-3}
m' \leq (m-1) (\tau(n,r)+1).
\end{equation}
By contradiction assume that $m' \geq (m-1) (\tau(n,r)+1)+1$. Let
$$
\tang(j)=\{l|1 \leq l \leq m', S_r(n, A_l) \;\;\ \mbox{\rm
  tangent to} \;\; S_r(n, A_j)\}.
$$
According to the definition of $\tau(n, r)$, we know the
cardinality of $\tang(j)$ is less than $\tau(n,r)$. We first pick
$j_1$ from $\{0,1,\cdots, m'\}$, then pick $j_2$ from
$\{0,1,\cdots, m'\}-\tang(j_1)$. And we continue this process by
always picking $j_{k+1}$ from
$$
\{0,1,\cdots, m'\}-\bigcup_{l=1}^k \tang(j_l).
$$
Since the cardinality of the above set is strictly greater than
$0$ when $k \leq m-1$, we can pick $j_1, j_2, \cdots, j_m$ from
the index set $\{1, 2, \cdots, m'\}$ such that $S_r(n,
A_{j_1}),S_r(n, A_{j_2}),\cdots,S_r(n, A_{j_m})$ are
non-overlapping and every two of them are not tangent to each
other. Then we can find a small enough real number $\varepsilon >
0$ and increase the radius $r$ to $r+\varepsilon$ such that
$$
S_{r+\varepsilon}(n, A_{j_1}),S_{r+\varepsilon}(n,
A_{j_2}),\cdots,S_{r+\varepsilon}(n, A_{j_m})
$$
are still non-overlapping. However this contradicts the
maximality of $r=r_{opt}(n, m)$.

The combination of the three formulas\eqr{formula-1},
\eqr{formula-2}, \eqr{formula-3} will lead to
\begin{multline}
  \frac{\int\!\!\!\int_{D_1} \prod_{j<k}{|e^{i \theta_j}-e^{i
        \theta_k}|}^2 d\theta_1 d\theta_2 \cdots
    d\theta_n}{\int\!\!\!\int_{D_1 \cap D_3} \prod_{j<k}{|e^{i
        \theta_j}-e^{i \theta_k}|}^2 d\theta_1 d\theta_2 \cdots
    d\theta_n}\\
  \leq \left(\frac{\int\!\!\!\int_{D_1} \prod_{j<k}{|e^{i
          \theta_j}-e^{i \theta_k}|}^2 d\theta_1 d\theta_2 \cdots
      d\theta_n}{\int\!\!\!\int_{D_1 \cap \tilde{D}_2} \prod_{j<k}{|e^{i
          \theta_j}-e^{i \theta_k}|}^2 d\theta_1 d\theta_2 \cdots
      d\theta_n}-1\right)(\tau(n,r)+1).
\end{multline}
Note that the inequality above is in fact stronger than the claim
in the theorem. We can reach the conclusion of the theorem by
relaxing the right hand side of the inequality (by ignoring
$-1$).
\end{proof}

\begin{co}
  When $m \rightarrow \infty$, asymptotically we have
  $$
  \Delta  (n, m) \geq 2\sqrt{n} r_0^E(n, m) \frac{1}{2}
  (\tau(2n^2-1)+1)^{-1/n^2}.
  $$

\end{co}

\begin{proof}
  We only sketch the idea of the proof. Intuitively $U_r^E(n,
  A_0)$ looks more ``flat'' when $m \rightarrow \infty$
  (consequently $r \rightarrow 0$), so $V(U_r^E(n, A_0))$ can be
  approximated by the volume of $U_r^E(n, A_0)$'s projection to
  the tangent space of $U(n)$ at $A_0$:
  $$
  \int\!\!\!\int_{D_1 \cap \tilde{D}_2} \prod_{j<k}{|e^{i
      \theta_j}-e^{i \theta_k}|}^2 d\theta_1 d\theta_2 \cdots
  d\theta_n \sim C (r_0^E)^{n^2}
  $$
  for some constant $C$. The same argument will lead to
  $$
  \int\!\!\!\int_{D_1 \cap D_3} \prod_{j<k}{|e^{i
      \theta_j}-e^{i \theta_k}|}^2 d\theta_1 d\theta_2 \cdots
  d\theta_n \sim C (2r_{opt})^{n^2}
  $$
  for the same constant $C$. For any fixed $n$, $\tau(n, r)$ will
  approach to the standard kissing number in Euclidean space
  $\tau(2n^2-1)$ when $r$ goes to zero. Combining the three
  approximations, we reach the claim according to the previous
  theorem.

\end{proof}

$U(n)$ is a compact Lie group equipped with a Riemannian metric.
Given two points $A_0, A_1 \in U(n)$, one can always find a
geodesic $\gamma(t)$ (mapping from $[0,1]$ to $U(n)$) which will
connect these two points, i.e. $\gamma(0)=A_0$ and
$\gamma(1)=A_1$. Recall that the Euclidean distance of $A_0$ and
$A_1$ is defined to be $\|A_0-A_1\|$. We further define the
Riemannian distance between $A_0$ and $A_1$ to be:
$$
\dist(A_0, A_1)=\int_0^1 \|\gamma'(t)\| dt.
$$
As a Lie group $U(n)$ is homogeneous. In particular one has
that
$$
\dist(A_0, A_1)=\dist(UA_0, UA_1)=\dist(A_0U, A_1U)
$$
for any $U \in U(n)$. The following theorem utilizes the
homogeneity and the relationship between the Riemannian distance
and Euclidean distance to derive another upper bound for the
diversity sum in general and it is the base of the second
approach.

\begin{thm} \label{Second-Approach}
  Let $f(\cdot)$ and $g(\cdot)$ be two fixed monotone increasing
  real functions. If  $$g(\|A_0-A_1\|) \leq \dist(A_0, A_1) \leq
  f(\|A_0-A_1\|)$$ for any two unitary matrices $A_0$ and $A_1$,
  then
  $$
  \Delta  (n, m) \leq g^{-1}(2f(r_0^E(n, m)))/(2\sqrt{n}).
  $$
\end{thm}
\begin{proof}
  For a fixed unitary constellation $\V=\{A_1, A_2, \cdots,
  A_m\}$, consider
  $$
  U_r^E(n, A_1), U_r^E(n, A_2), \cdots, U_r^E(n, A_m)
  $$
  for $r>0$. We can increase $r$ until there exist $j$ and
  $k$ such that $U_r^E(n, A_j)$ and $U_r^E(n, A_k)$ are tangent to each other at a point $A_0$.
  As examined in Theorem~\ref{Sphere-Upper-Bound}, one can make a
  conclusion that $r \leq r_0^E(n,m)$. Accordingly we have
  \begin{multline*}
    \dist(A_j, A_k) \leq \dist(A_j, A_0)+\dist(A_k, A_0) \\
    \leq f(\|A_j-A_0\|)+f(\|A_k-A_0\|)=2f(r) \leq 2f(r_0^E(n,
    m)).
  \end{multline*}
  On the other hand since $g$ is monotonically increasing one
  has:
  $$
  \|A_j-A_k\| \leq g^{-1}(\dist(A_j, A_k)).
  $$
  The combination of the above two inequalities will lead to
  $$
  \|A_j-A_k\| \leq g^{-1}(2f(r_0^E(n, m))).
  $$
  Immediately we will have
  $$
  \sum \V \leq g^{-1}(2f(r_0^E(n, m)))/(2\sqrt{n}).
  $$
  Since $\V$ is an arbitrary unitary constellation, the claim
  in the theorem follows.
\end{proof}

Based on the above theorem, the following corollary gives the
second upper bound (B2).
\begin{co}
  For a real number $r$, let $\lfloor r \rfloor$ denote the
  greatest integer less than or equal to~$r$, then
\begin{equation}
\Delta  (n, m) \leq \sin \sqrt{\frac{\pi^2}{n}\left\lfloor
\frac{(r_0^E)^2 (n, m)}{4}
\right\rfloor +\frac{4}{n}
\arcsin^2\sqrt{\frac{(r_0^E)^2(n, m)}{4}-\left\lfloor
\frac{(r_0^E)^2(n, m)}{4} \right\rfloor}}. \tag{B2}
\end{equation}
\end{co}

\begin{proof}
  Consider $I$ and another point $U = V \diag( e^{i \theta_1},
  e^{i \theta_2}, \cdots, e^{i \theta_n}) V^* $, where $ -\pi \leq
  \theta_j < \pi$. It is known that~\cite{ed99} the geodesic from
  $I$ to $U$ can be parameterized by
  $$
  \gamma(t)=V \diag(e^{i\theta_1 t},e^{i\theta_2 t}, \cdots,
  e^{i\theta_n t}) V^*,
  $$
  where $0 \leq t \leq 1$. The Riemannian distance from $I$ to
  $U$ is
  $$
  \dist(I,U)=\sqrt{\theta_1^2+\theta_2^2+\cdots+\theta_n^2}.
  $$
  We want to derive $g(\cdot), f(\cdot)$ as in
  Theorem~\ref{Second-Approach}. Suppose the Euclidean distance
  between $I$ and $U$ is $r$, i.e.,
  $$
  \sin^2 \frac{\theta_1}{2}+\sin^2
  \frac{\theta_2}{2}+\cdots+\sin^2 \frac{\theta_n}{2}= r^2/4.
  $$
  After substituting with $x_j=\sin^2 \theta_j/2$ and denoting
  $G(x)=\arcsin^2 \sqrt{x}$, we convert the above problem to the
  following optimization problem:

  Find the minimum and maximum of the function
  $$
  F(x_1, x_2, \cdots,
  x_n)=\theta_1^2+\theta_2^2+\cdots+\theta_n^2=4(G(x_1)+G(x_2)+\cdots+G(x_n))
  $$
  with the constraints $x_1+x_2+\cdots+x_n=r^2/4$ and $0 \leq
  x_j \leq 1$ for $j=1, 2, \cdots, n$. Since $G(x)$ is a convex
  function on $[0,1]$, we derive the lower bound of $F(x_1, x_2,
  \cdots, x_n)$,
\begin{equation}  \label{Riemannian-Euclidean}
4 n \arcsin^2 (r/(2 \sqrt{n})) \leq F(x_1, x_2, \cdots, x_n).
\end{equation}
In the sequel we are going to calculate the upper bound of
$F(x_1, x_2, \cdots, x_n)$. Without loss of generality, we assume
$0 \leq x_1 \leq x_2 \leq \cdots \leq x_n \leq 1$. Let $k=
\lfloor r^2/4 \rfloor$ and $\alpha= r^2/4 -k$, we claim that
$F(x_1, x_2, \cdots, x_n)$ will reach its maximum when
$$
x_j=\left \{ \begin{tabular}{cc}
  $0$ & $1 \leq j \leq n-k-1$ \\
  $\alpha$ & $j=n-k$ \\
  $1$ & $n-k+1 \leq j \leq n$ \\
\end{tabular} \right.
$$
Suppose by contradiction that $F$ reaches its maximum at
$(x_1, x_2, \cdots, x_n)$ with $x_1 > 0$. Now from
$$
x_1+x_{n-k}+x_{n-k+1}+\cdots+x_n \leq r^2/4=k+\alpha,
$$
surely one can find $x'_{n-k}, x'_{n-k+1}, \cdots, x'_n$ such
that
$$
x_1+x_{n-k}+x_{n-k+1}+\cdots+x_n=x'_{n-k}+x'_{n-k+1}+\cdots+x'_n,
$$
with $x'_j \geq x_j$ for $j=n-k, n-k+1, \cdots, n$. Now set
$x^*_1=0$, $x^*_j=x_j$ for $j=2, 3, \cdots, n-k-1$ and
$x^*_j=x'_j$ for $j=n-k, n-k+1, \cdots, n$. By the mean value
theorem, there exist $\zeta_j$'s with $x^*_1=0 \leq \zeta_1 \leq
x_1$ and $x_j \leq \zeta_j \leq x^*_j$ for $j=2, 3, \cdots, n$
such that
$$
F(x^*_1, x^*_2, \cdots, x^*_n)-F(x_1, x_2, \cdots,
x_n)=\sum_{j=1}^n G'(\zeta_j) (x^*_j-x_j).
$$
Since $G(x)$ is a strictly convex function, we have
$$
0 < G'(\zeta_1) < G'(\zeta_2) < \cdots < G'(\zeta_n).
$$
Now
\begin{multline*}
  F(x^*_1, x^*_2, \cdots, x^*_n)-F(x_1, x_2, \cdots, x_n) \geq
  G'(\zeta_2)(\sum_{j=2}^n (x^*_j-x_j))
  -G'(\zeta_1)(x_1-x^*_1)\\
  =(G'(\zeta_2)-G'(\zeta_1))(x_1-x^*_1)
  =(G'(\zeta_2)-G'(\zeta_1))x_1>0.
\end{multline*}
This contradicts the maximality of $F$ at $(x_1,x_2, \cdots,
x_n)$.  Applying exactly the same analysis to $x_2, x_3, \cdots,
x_{n-k-1}, x_{n-k}$ we deduce that $x_j=0$ for $j=2, 3, \cdots,
n-k-1$ and $x_{n-k}=\alpha$. So the upper bound of $F$ can be
given as
$$
F(x_1, x_2, \cdots, x_n) \leq 4\left( k
  \frac{\pi^2}{4}+\arcsin^2(\sqrt{\alpha})\right).
$$
Take $g(r)=2 \sqrt{n} \arcsin (r/(2 \sqrt{n}))$ and $f(r)=2
\sqrt{k \pi^2/4+\arcsin^2\sqrt{\alpha}}$, the corollary follows
according to the previous theorem.
\end{proof}

\begin{figure}
\begin{center}
  \includegraphics[totalheight=2.3in]{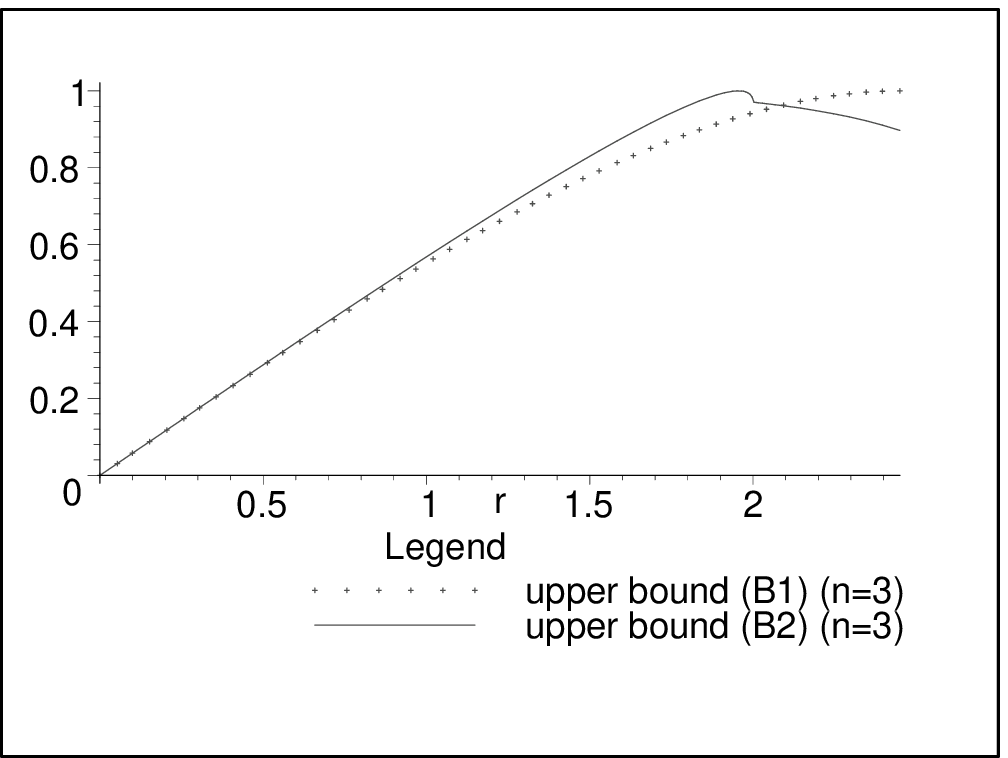}
  \includegraphics[totalheight=2.3in]{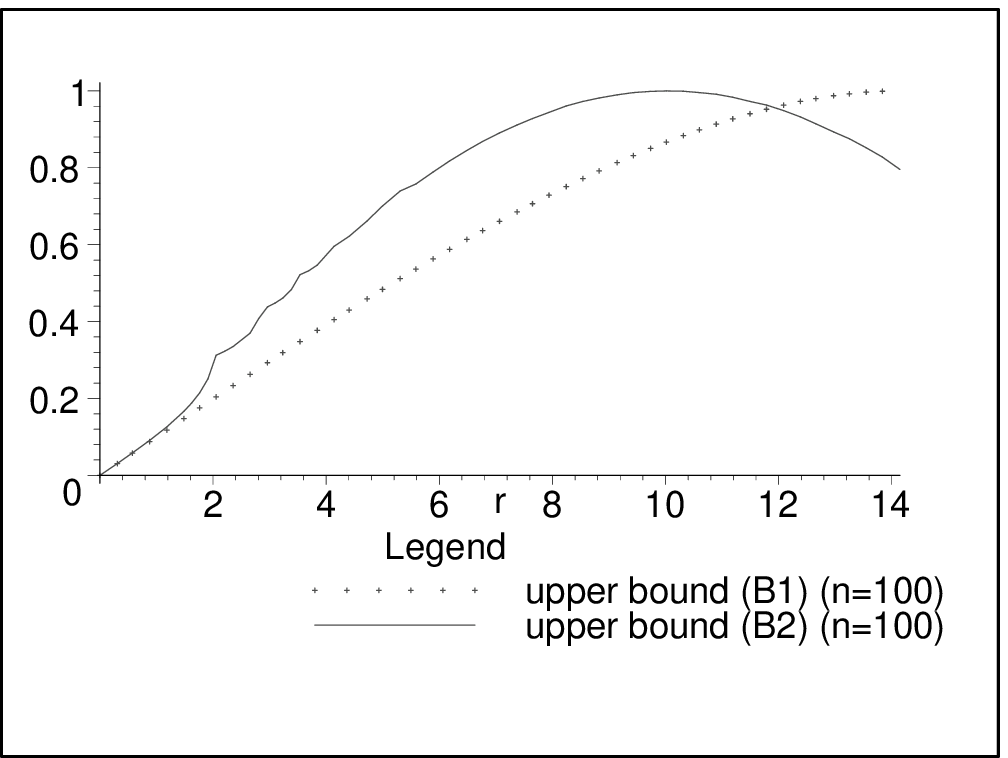}
\caption{The comparisons of two upper bounds as functions for
  $n=3$ and $n=100$} \label{two-upper-bound}
\end{center}
\end{figure}

Note that both upper bound (B1) and upper bound (B2) depend on
$r_0^E(n,m)$. In Figure~\ref{two-upper-bound} we plot both upper
bounds as functions of $r_0^E(n, m)$ for $3$ and $100$
dimensions.  One can see that if and only if $r_0^E(3,m)>2.0881$,
the upper bound (B2) is tighter than the upper bound (B1). While
for the $100$ dimension case, the upper bound (B1) is tighter
than the upper bound (B2) if and only if $r_0^E(100, m)>11.9155$.
In fact it can be checked that asymptotically when $n$ is large
enough, upper bound (B2) is tighter than upper bound (B1) if and
only if $r_0^E(n, m)> 1.1892 \sqrt{n}$.

For a packing problem on a manifold, alternatively one can choose
the neighborhood to be a small ``ball'' with Riemannian radius
$r$. This will be our third approach to derive an upper bound for
the diversity sum. For a particular $A \in U(n)$, let
$$
U_r^R(n, A)=\{U \in U(n)|\dist(U, A) \leq r\}.
$$
Note that the constraint $\dist(U,I) \leq r$ is equivalent to
$$
\theta_1^2+\theta_2^2+\cdots+\theta_n^2 \leq r^2.
$$
Therefore we apply the same argument as in the proof of
Theorem~\ref{Euclidean-Ball} and conclude that:
$$
V(U_r^R(n))=\frac{\int\!\!\!\int_{D_1 \cap D_4}
  \prod_{j<k}{|e^{i \theta_j}-e^{i \theta_k}|}^2 d\theta_1
  d\theta_2 \cdots d\theta_n }{\int\!\!\!\int_{D_1}
  \prod_{j<k}{|e^{i \theta_j}-e^{i \theta_k}|}^2 d\theta_1
  d\theta_2 \cdots d\theta_n} V(U(n)),
$$
where $D_1$ was defined in\eqr{D1} and
\begin{equation}                           \label{D4}
  D_4:=\{(\theta_1, \theta_2, \cdots, \theta_n)| \sum_{j=1}^n
\theta_j^2 \leq r^2 \}.
\end{equation}
Instead of considering the Euclidean neighborhoods $U_r^E(n,
A_1), U_r^E(n, A_2), \cdots, U_r^E(n, A_m)$, we can consider the
Riemannian neighborhood $U_r^R(n, A_1), U_r^R(n, A_2), \cdots,
U_r^R(n, A_m)$. Utilizing the fact that the Euclidean distance
$\|A_j-A_k\|$ and the Riemannian distance $\dist(A_j, A_k)$ are
related (compare with Formula\eqr{Riemannian-Euclidean}):
$$
4 n \arcsin^2 (\|A_j-A_k\|/(2 \sqrt{n})) \leq \dist(A_j, A_k)
$$
for any two unitary matrices $A_j$ and $A_k$, we can derive
the third upper bound (B3). The proof of the following theorem is
very similar to the one of Theorem~\ref{Second-Approach} and for
the sake of brevity we omit it.
\begin{thm}
  Let $D_1$ and $D_4$ be defined as in\eqr{D1} and \eqr{D4} and
  assume $ r_0^R(n, m) $ is the solution to the following
  equation (with variable $r$):
\begin{equation} \label{Riemannian-Equality}
m\int\!\!\!\int_{D_1 \cap D_4} \prod_{j<k}{|e^{i \theta_j}-e^{i
\theta_k}|}^2 d\theta_1 d\theta_2 \cdots d\theta_n =
\int\!\!\!\int_{D_1} \prod_{j<k}{|e^{i \theta_j}-e^{i
\theta_k}|}^2 d\theta_1 d\theta_2 \cdots d\theta_n,
\end{equation}
then
\begin{equation}
\Delta  (n, m) \leq \sin \left(\frac{r_0^R(n, m)}{\sqrt{n}}\right). \tag{B3}
\end{equation}
\end{thm}

\begin{figure}[ht]
  \centerline{\psfig{figure=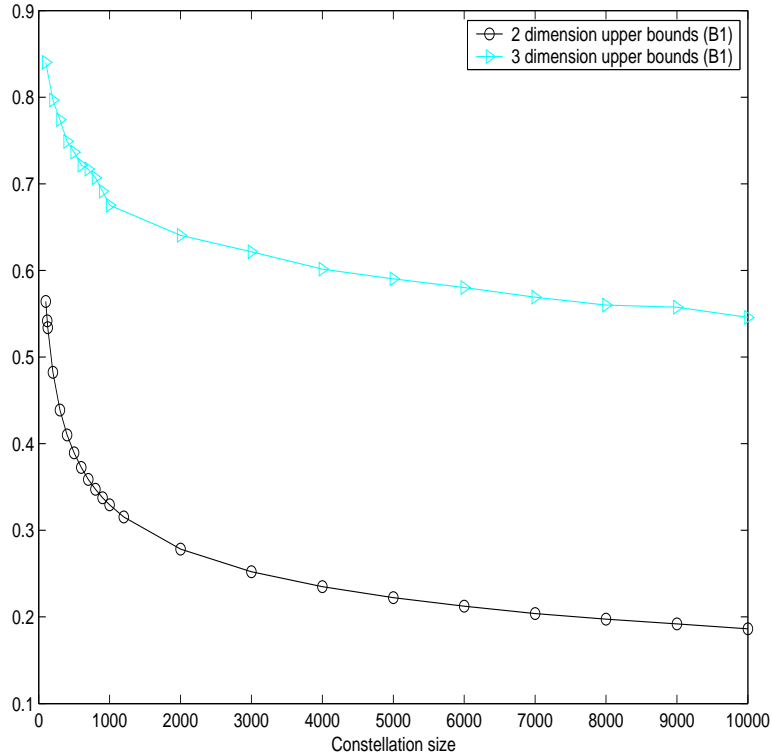,width=4in,height=4in}}
\caption{Upper bounds for 2 and 3 dimensional constellations}
\label{Two-Bound}
\end{figure}

We gave three approaches to derive upper bounds for the diversity
sum and hence also for the diversity product.  All of them
involve the calculation of $r_0^E(n, m)$ or $r_0^R(n, m)$, which
are the solutions of equation\eqr{Euclidean-Equality} and
equation\eqr{Riemannian-Equality}, respectively. Fortunately we
are dealing with finding a root of a monotone increasing function
(recall that both $V(U_r^E(n, m))$ and $V(U_r^R(n, m))$ are
monotone increasing functions with respect to $r$), the bisection
method~\cite{at78b} will be highly effective to solve this kind
of problem. Our numerical experiments for small size
constellations with small dimensions show that upper bound (B3)
is looser than the first two upper bounds. However when $m$ goes
to infinity, these three upper bounds give almost the same
estimation. This makes sense because asymptotically the small
balls look like a $n^2$ dimensional ball in Euclidean space. One
can see the derived upper bounds for $2$ and $3$ dimensional
constellations in Figure~\ref{Two-Bound}.

We compare the derived upper bounds with the currently existing
one presented in~\cite{li02}. For $n=2$ the upper bounds derived
by Liang and Xia~\cite{li02} tend to be better when $m \leq 100$
and our bounds become tighter when $m \geq 100$ (see the
following Table~\ref{UB-table}). For $n \geq 3$ Liang and
Xia~\cite{li02} outlined a method by considering a sphere packing
computation in $S(n)$. It is our belief that this method will
result in a weaker bound than the upper bounds we derived in this
paper. For the sample programs to do the upper bound calculation,
we refer to~\cite{ha03u2}.

\vspace{0.3cm}
\begin{ta}  \label{UB-table}
  For $n=2$ the following table compares the upper bounds
  in~\cite{li02} with our new bounds (B1) and (B2).\\
\begin{tabular}{|c|c|c|c|c|c|c|c|c|}
  \hline
  m&24 &48 &64 &80 &100 &120 &128&1000 \\ \hline
  \!\!upper bounds in~\cite{li02}\!\!&0.6746&0.6193 &0.5969
  &0.5799 &0.5632 &0.5499 &0.5452& \\ \hline
  upper bound (B1) &0.7598&0.6603 &0.6131 &0.5932 &0.5578
  &0.5425 &0.5347&0.3270 \\ \hline
  upper bound (B2) &0.7794&0.6734 &0.6235 &0.6026 &0.5654
  &0.5496 &0.5415&0.3285\\
  \hline
\end{tabular}
\end{ta}

\vspace{0.3cm}

One interesting fact about the limiting behavior of $\Delta (n,
m)$ (when $m \rightarrow \infty$) is its connection to the Kepler
problem~\cite{co93b}. Certainly one can use Kepler
density~\cite{co93b} to obtain a tighter bound of the diversity
sum asymptotically.

\Section{Conclusions and Future Work}

We presented three approaches to derive upper bounds for the
diversity sum of unitary constellations of any dimension $n$ and
any size $m$. The derived bounds seem to improve the existing
bounds when $n=2$ and $m\geq 100$. When $n$ is large the exact
computation of $r_0^E$ is rather involved and hence it is also
computationally difficult to compute the bounds (B1) and (B2).
Nonetheless it is our belief that the resulting upper bounds (B1)
and (B2) become fairly tight as soon as $m$ is sufficiently
large.

It was pointed out that the resulted upper bounds also apply for
the diversity product, although the bounds seem to be less tight
in this situation.  The future work may involve the derivation of
a tighter upper bound analysis for the diversity product of
unitary constellations using differential geometric means.

\subsection*{Acknowledgments}

We started this research while spending a month at the Institute
Mittag-Leffler in Stockholm in May 2003. The hospitality and the
financial support of the Institute Mittag-Leffler are greatly
acknowledged. We are also grateful for stimulating discussions we
had on this subject with Professor Uwe Helmke and Professor
Xuebin Liang.


\end{document}